\documentclass[12pt,twoside]{amsart}
\usepackage{amssymb}

\nonstopmode

\numberwithin{equation}{section}

\font\tengothic=eufm10 scaled\magstep 1
\font\sevengothic=eufm7 scaled\magstep 1
\newfam\gothicfam
      \textfont\gothicfam=\tengothic
      \scriptfont\gothicfam=\sevengothic

\newtheorem{theorem}{Theorem}[section]
\newtheorem{lemma}[theorem]{Lemma}
\newtheorem{proposition}[theorem]{Proposition}
\newtheorem{corollary}[theorem]{Corollary}

\theoremstyle{definition}

\newtheorem{definition}[theorem]{Definition} 
\newtheorem{remark}[theorem]{Remark}
\newtheorem{example}[theorem]{Example}

\newtheorem{notation}[theorem]{Notation}

\newtheorem{notation and remark}[theorem]{Notation and Remark}

\newcommand\rank{\operatorname{rank}}

\newcommand{\kindalong}{\ \hbox{$\hbox to .35in{\rightarrowfill}$} \ }
\newcommand{\reallylong}{\hbox{$\hbox to .5in{\rightarrowfill}$}  }

\def\ra{\rightarrow}

\def\End{{\rm End}}
\def\In{{\rm In}}

\def\ba{{\bf a}}
\def\bff{{\bf f}}
\def\Min{{\rm Min}}
\def\ra{\longrightarrow}

\def\deg{{\rm deg }}
\def\End{{\rm End}}

\def\rank{{\rm rank}}

\def\ol{\overline}
\def\ti{\tilde}
\def\ni{\noindent}

\def\Gr{{\rm Gr}}

\def\In{{\rm In}}

\def\Min{{\rm Min}}
\def\ra{\longrightarrow}

\def\deg{{\rm deg }}
\def\End{{\rm End}}

\def\rank{{\rm rank}}

\def\Max{{\rm Max}}
\def\bff{{\rm f}}

\def\Gr{{\rm Gr}}
\def\bff{{\bf f}}

\newfont{\bg}{cmr10 scaled\magstep4}
\newcommand{\bigzerol}{\smash{\hbox{\bg 0}}}
\newcommand{\bigzerou}{\smash{\lower1.7ex\hbox{\bg 0}}}
\newcommand{\bigast}{\smash{\lower3ex\hbox{\bg *}}}

\begin{document}


\begin{center}
{\Large \bf The central simple modules of \\[0.5ex] Artinian Gorenstein algebras}
\end{center}
\medskip\medskip\medskip

\begin{center}
{\sc Tadahito Harima} \\[1ex]
Department of Mathematics, Hokkaido University of Education, \\
Kushiro  085-8580, Japan \\
E-mail: harima@kus.hokkyodai.ac.jp 
\medskip\medskip

{\sc Junzo Watanabe} \\[1ex]
Department of Mathematics, Tokai University, \\
Hiratsuka 259-1292, Japan \\
E-mail: junzowat@keyaki.cc.u-tokai.ac.jp
\end{center}
\medskip\medskip\medskip


\begin{quote}
{\footnotesize 
{\sc Abstract.}\ 
Let $A$ be a standard graded Artinian $K$-algebra, with char $K=0$. 
We prove the following. 
\begin{itemize}
\item[1.] 
$A$ has the Weak Lefschetz Property (resp. Strong Lefschetz Property) 
if and only if $\Gr_{(z)}(A)$ has the Weak Lefschetz Property 
(resp. Strong Lefschetz Property) for some linear form $z$ of $A$. 
\item[2.] 
If $A$ is Gorenstein, 
then $A$ has the Strong Lefschetz Property
if and only if there exists a linear form $z$ of $A$ such that 
all central simple modules of $(A,z)$ have the Strong Lefschetz Property. 
\end{itemize}
As an application of these theorems,  
we give some new classes of Artinian complete intersections  
with the Strong Lefschetz Property. 
\medskip\medskip

}
\end{quote}
\medskip\medskip\medskip\medskip


\section{Introduction}
Let $A=\oplus_{i=0}^c A_i$ be a standard graded Artinian $K$-algebra, 
where $A_c \neq (0)$. 
An algebra $A$ has the {\em Weak Lefschetz Property} (WLP) 
if there exists a linear form $g \in A_1$ such that 
the multiplication  
$$
\times g: A_i \ra A_{i+1}
$$
is either injective or surjective for every $i=0,1,\ldots,c-1$. 
Such an element $g$ will be called a {\em weak Lefschetz element}.  
An algebra $A$ has the {\em Strong Lefschetz Property} (SLP) 
if there exists a linear form $g \in A_1$ such that 
the multiplication  
$$
\times g^{c-2i}: A_i \ra A_{c-i}
$$
is bijective for all $i=0,1,\ldots,[c/2]$. 
In this case, such an element $g$ is called a {\em strong Lefschetz element}. 
The set of weak or strong Lefschetz elements $g\in A_1$ is 
a Zariski open subset of $A_1$. 

These are fundamental properties 
and have been studied by many authors. 
Below are some of the results so far obtained.   
\begin{itemize}
\item[(1)]
The SLP is preserved by tensor product (\cite{jW87a}). 
\item[(2)]
If $A$ has the SLP and if a homogeneous element $f$ is general enough, 
then $A/0:f$ has the SLP (\cite{jW87a}).
\item[(3)] 
Let $R$ be a polynomial ring and $I$ a homogeneous ideal. 
Then $R/I$ has the WLP (resp. SLP) 
if and only if $R/{\rm Gin}(I)$ has the WLP (resp. SLP), 
where ${\rm Gin}(I)$ is the generic initial ideal of $I$ 
with respect to the reverse lexicographic order (\cite{aW03}).
\item[(4)]
If $R/{\rm In}(I)$ has the WLP (resp. SLP), 
then $R/I$ has the WLP (resp. SLP), 
where ${\rm In}(I)$ is the initial ideal of $I$ with respect to a term order (\cite{aW03}).
\item[(5)]
The finite free extension of an Artinian $K$-algebra 
with the SLP has again the SLP (\cite{tHjW03}, \cite{tHjW04}), 
if the fiber has the Strong Lefschez Property. 
This is a generalization of (1). 
\end{itemize}

These results have widened the class of Artinian $K$-algebras which we know  have 
the WLP or SLP. 
In this paper we would like to prove two theorems in this direction. 
First we prove the following.  

\begin{theorem} \label{main-th1}   
Let $K$ be a field of characteristic zero 
and let $A$ be a standard graded Artinian $K$-algebra. 
Then 
$A$ has the WLP (resp. SLP) if and only if 
$\Gr_{(z)}(A)$ has the WLP (resp. SLP) for some linear form $z\in A$. 
\end{theorem}
\noindent
The ``if" is the essential part. 
The key is to show that 
if $g_1$ and $g_2$ are general linear forms of $A$ and $\Gr_{(z)}(A)$ respectively, 
then $$\dim_K A/g_1A\leq\dim_K\Gr_{(z)}(A)/g_2\Gr_{(z)}(A).$$
This can be proved, somewhat surprisingly, by looking at the Artinian algebras 
as modules over a two dimensional polynomial ring. 
(See Proposition~\ref{prop}.)   

Using Theorem~\ref{main-th1}, 
we prove the main theorem of this paper as follows.  

\begin{theorem} \label{main-th2}   
Let $K$ be a field of characteristic zero 
and let $A$ be a standard graded Artinian Gorenstein $K$-algebra. 
Then the following conditions are equivalent.  
\begin{itemize}
\item[(i)] 
$A$ has the SLP. 
\item[(ii)]
There exists a linear form $z$ of $A$ 
such that all the central simple modules of $(A,z)$ have the SLP. 
\end{itemize}
\end{theorem}

The definition of the ``central simple module'' is given 
in Section~4.  These are defined for a pair $(A, z)$ of Artinian 
$K$-algebra $A$ and a linear form $z \in A$. 
 
Let 
\[
\times : A \to \End(A)
\]
be the regular representation of $A$, so that 
$\times a$, for $a \in A$, is the endomorphism of $A$ defined by $(\times a) x= ax$. 
Put $J=\times z$ and let 
$$C(J)=\{X \in \End(A) | XJ=JX   \}$$ 
be the centralizer of $J$ in $\End(A)$. 
It is easy to determine the algebra $C(J)$ as a set in the matrix ring $\End(A)$,  
since the condition $XJ=JX$ is a system of linear equations in the entries of $X$ as 
unknowns. 
Some results on $C(J)$ are found in 
\cite{rB00}, \cite{fG59}, and \cite{hTaA32}. 
It turns out  that
the central simple modules of $(A,z)$ are simple modules of
$C(J)$ and moreover they exhaust all the isomorphism types of
such modules. (For details see \cite{tHjW0000}.)

In the present paper one does not consider the centralizers.  Nonetheless
the knowledge of $C(J)$ might be helpful to understand the proof of  
Theorem~\ref{main-th2}.  
A systematic approach to the central simple 
modules is made in \cite{tHjW0000}. 
In Section 4 it will be shown that all central simple modules for 
Gorenstein $K$-algebras have symmetric Hilbert functions (Proposition~\ref{important_property}). 
This was really the motivation for considering these modules. 
It enables us to make an inductive argument when we try to prove certain types of 
Gorenstein $K$-algebras have the SLP.  
Some such examples will be given in Section~6. 

It should be noted that the ``SLP'' is defined not only for 
algebras but also for modules.  But in any case we assume the 
symmetry of the Hilbert function for those modules with the  
SLP.  The important fact is that then the SLP is preserved by the 
tensor product over a field of characteristic zero (\cite{jW87a}, Corollary 3.5). 
It should be noted, moreover,  that even the WLP often fails to 
hold without this restriction. 
Thus all propositions concerning WLP as well as SLP are stated over a field  $K$  of 
characteristic zero.

%
%
%
%

\section{Some basic results on the WLP and SLP}

Let $A=\oplus_{i=0}^cA_i$ be a standard graded Artinian $K$-algebra. 
We define 
the {\em Sperner number} of $A$ by 
$$
{\rm Sperner}(A)={\rm Max}\{\dim A_i\}
$$
and {\em CoSperner number} of $A$ by 
$$
{\rm CoSperner}(A)=\sum _{i=0} ^{c-1} {\rm Min}\{\dim A_i, \dim A_{i+1}\}.
$$
If $A$ has the WLP then the Hilbert function of $A$ is unimodal, and  
we have 
$$
{\rm CoSperner}(A)=\dim A-{\rm Sperner}(A). 
$$
Furthermore, if $A$ has the SLP, 
the Hilbert function of $A$ is an SI-sequence. 
(For details see \cite{tHjMuNjW01}. 
In fact, in this paper,  
we gave a complete characterization of the Hilbert functions that 
can occur for Artinian $K$-algebras with the WLP or SLP.) 
For a homogeneous form $f \in A$, 
we denote the rank of the multiplication $\times f: A \ra A$ by $\rank(\times f)$. 
\medskip

The following is a characterization of the WLP. 

\begin{lemma}[\cite{jW87a}, Proposition 3.2] \label{lemma0}  
Let $A=\oplus_{i=0}^cA_i$ be a standard graded Artinian $K$-algebra, where $A_c\neq (0)$. 
Then we have the following. 
\begin{enumerate}
\item[$(1)$]
$\dim A/\ell A \geq {\rm Sperner}(A)$ 
\ for any linear form $\ell \in A_1$. 
\item[$(2)$]
$\rank (\times \ell)\leq{\rm CoSperner}(A)$ 
\ for any linear form $\ell \in A_1$. 
\item[$(3)$]
$A$ has the WLP if and only if  
$\dim A/gA = {\rm Sperner}(A)$ 
\ for a linear form $g \in A_1$. 
\item[$(4)$]
A has the WLP if and only if 
$\rank (\times g)={\rm CoSperner}(A)$ 
\ for a linear form $g \in A_1$. 
\end{enumerate}
\end{lemma}

Let $A=\oplus_{i=0}^c A_i$ be a standard graded Artinian $K$-algebra 
such that the Hilbert function of $A$ is an SI-sequence. 
Then, for $1\leq k\leq c$, put 
$$
{\rm SP_k}(A)=\sum_{i=0}^c\Max\{\dim A_i - \dim A_{i-k}, \ 0 \}, 
$$
where $\dim_K A_j=0$ for all $j<0$. 
We call 
$$
{\bf SP}(A)=({\rm SP_1}(A), {\rm SP_2}(A),\ldots,{\rm SP_c}(A))
$$
the {\em Sperner vector} of $A$. 
Note that ${\rm SP_1}(A)$ is equal to the Sperner number of $A$. 
\medskip
 
Using this we can characterize the SLP as follows. 

\begin{lemma} \label{lemma2}  
Let $A=\oplus_{i=0}^c A_i$ be a standard graded Artinian $K$-algebra, 
where $A_c\neq (0)$. 
Assume that the Hilbert function of $A$ is an SI-sequence. 
Then we have the following. 
\begin{enumerate}
\item[$(1)$]
$\dim_K A/f A \geq {\rm SP_k}(A)$ 
\ for any homogeneous form $f$ of degree $k$, 
where $1\leq k\leq c$. 
\item[$(2)$]
$\rank (\times f)\leq\dim A-{\rm SP_k}(A)$ 
\ for any homogeneous form $f$ of degree $k$, 
where $1\leq k\leq c$. 
\item[$(3)$]
$A$ has the SLP if and only if 
there exists a linear form $g\in A_1$ such that 
$$
\dim A/g^k A = {\rm SP_k}(A)
$$ 
for all $k=1,2,\ldots,c$. 
\item[$(4)$]
A has the SLP if and only if 
there exists a linear form $g\in A_1$ such that 
$$
\rank (\times g^k)=\dim A-{\rm SP_k}(A)
$$ 
for all $k=1,2,\ldots,c$. 
\end{enumerate}
\end{lemma} 

\begin{proof}
For all integers $k=1,2,\ldots,c$, 
put $u_k={\rm Max}\{i\mid \dim A_i\geq\dim A_{i-k}\}$. 
Then, noting that the Hilbert function of $A$ is unimodal and symmetric, 
it follows that 
$$
{\rm SP_k}(A) = \dim A_0+\cdots+\dim A_{k-1}+\sum_{i=k}^{u_k}(\dim A_i - \dim A_{i-k}). 
$$
Hence we have that
$$
\begin{array}{rl}
     & \dim A/fA \\[1ex]
=    & \dim A_0+\cdots+\dim A_{k-1}+
\sum_{i=k}^{u_k} \dim A_i/fA_{i-k} + 
\sum_{i=u_k+1}^c \dim A_i/fA_{i-k}  \\[1ex]
=    & \dim A_0+\cdots+\dim A_{k-1}+
\sum_{i=k}^{u_k}(\dim A_i-\dim fA_{i-k}) + 
\sum_{i=u_k+1}^c \dim A_i/fA_{i-k}  \\[1ex]
\geq & \dim A_0+\cdots+\dim A_{k-1}+
\sum_{i=k}^{u_k}(\dim A_i-\dim A_{i-k}) + 
\sum_{i=u_k+1}^c \dim A_i/fA_{i-k}  \\[1ex]
=    & {\rm SP_k}(A) + \sum_{i=u_k+1}^c \dim A_i/fA_{i-k}.
\end{array}
$$
Thus we obtain the inequality (1). 
In particular, 
$$
\begin{array}{rcl}
\mbox{$\dim A/fA={\rm SP_k}(A)$} & \Leftrightarrow & 
\times f:A_{i} \ra A_{i+k}\ \ \mbox{is} 
\ \left\{
\begin{array} {l}
\mbox{injective \ if}\ \dim A_i\leq \dim A_{i+k}, \\
\mbox{surjective \ if}\ \dim A_i\geq \dim A_{i+k}. 
\end{array}
\right.
\end{array}
$$
Furthermore, noting the equivalence of conditions 
$$
\begin{array}{rcl}
\mbox{$A$ has the SLP} & \Leftrightarrow & 
\exists g\in A_1\ \ \mbox{such that, for all $i < j$, } \\[1ex]
 & & \hspace{0.8cm}
\times g^{j-i}:A_{i} \ra A_{j}\ \ \mbox{is} 
\ \left\{
\begin{array}{l}
\mbox{injective \ if}\ \dim A_i\leq \dim A_{j}, \\
\mbox{surjective \ if}\ \dim A_i\geq \dim A_{j},  
\end{array}
\right.
\end{array}
$$
the assertion (3) easily follows from the argument above. 
The assertions (2) and (4) easily follow from (1) and (3), 
respectively.  
\end{proof}

\begin{lemma}\label{lemma3}   
Let $A$ and $B$ be two standard graded Artinian $K$-algebras 
with the same Hilbert function. 
Let $g$ and $g'$ be linear forms of $A$ and $B$, respectively. 
Assume the Jordan canonical form of $\times g: A\ra A$ is 
the same as that of $\times g':B\ra B$. 
Then $A$ has the WLP (resp. SLP) if and only if 
$B$ has the WLP(resp. SLP). 
\end{lemma}

\begin{proof}
This follows from Lemma \ref{lemma0} (4) and Lemma \ref{lemma2} (4).  
\end{proof}

\begin{lemma}[\cite{tHjW03}, Proposition 18]\label{lemma4}   
Let $A$ be a standard graded Artinian $K$-algebra 
with a symmetric Hilbert function, 
where ${\rm char}(K)=0$. 
Let $u$ be a new variable. 
Then the following conditions are equivalent. 
\begin{enumerate}
\item[$(1)$]
$A$ has the SLP. 
\item[$(2)$]
$A[u]/(u^\alpha)$ has the WLP for all positive integers $\alpha$. 
\end{enumerate}
\end{lemma}

%
%
%
%

\section{Proof of Theorem \ref{main-th1}}

We need some preparations for the proof of Theorem~\ref{main-th1}. 

\begin{notation and remark} \label{notation and remark 3-1}  
Let $A=\oplus_{i=0}^cA_i$ be a standard graded Artinian $K$-algebra. 
For any linear form $z \in A_1$, 
consider the associated graded ring 
$$
\Gr_{(z)}(A)=A/(z)\oplus (z)/(z^2) \oplus (z^2)/(z^3) \oplus \ldots \oplus (z^{p-1})/(z^p), 
$$
where $p$ is the least integer such that $z^p  =0$. 
As is well known $\Gr_{(z)}(A)$ is endowed with a commutative ring structure. 
The multiplication in $\Gr_{(z)}(A)$ is given by 
$$
(a+(z^{i+1}))(b+(z^{j+1}))=ab+(z^{i+j+1}), 
$$
where $a \in (z^i)$ and $b \in (z^j)$. 

For any homogeneous form $f \in R=K[x_1,\ldots,x_n]$,  
it is possible to write uniquely 
$$
f= f_0 + f_1x_n + f_2 x_n^2 + \cdots + f_k x_n^k, 
$$
where  $f_i$  is a homogeneous form in $K[x_1, \ldots, x_{n-1}]$.
Denote by ${\rm In}'(f)$ the term $f_jx_n^{j}$ 
for the minimal $j$ such that $f_j \not  = 0$.
Furthermore we define 
${\rm In}'(I)$ to be the homogeneous ideal of $R$ 
generated by the set $\{ {\rm In}' (f) \}$,  
where $f$ runs over homogeneous forms in $I$. 
Suppose that $n=2$.  
Then one notices that 
${\rm In}'(I)$ coincides with the initial ideal ${\rm In}(I)$ of $I$ 
with respect to the reverse lexicographic order with $x_1 > x _2$. 

Let $z$ be the image of $x_n$ in $A=R/I$. 
Then it is easy to show that ${\rm Gr}_{(z)}(A) \cong R/{\rm In}'(I)$. 
So we consider that $\Gr_{(z)}(A)$ inherits a standard grading from $R$. 
Let ${\rm In}(I)$ be the initial ideal of $I$ 
with respect to the reverse lexicographic order with $x_1 > \cdots > x_n$. 
Noting that ${\rm In}({\rm In}'(I))={\rm In}(I)$, 
we have that all of $\Gr_{(z)}(A)$, $R/{\rm In}(I)$ and $A=R/I$ 
have the same Hilbert function. 
\end{notation and remark}

\begin{notation and remark} \label{notation and remark 3-2}   
Let $S=K[Y,Z]$ be the polynomial ring over an infinite field $K$.  
Let $V$ be a graded Artinian $S$-module. 
Write $V\cong F/N$, 
where 
$$
F=\oplus_{j=1}^s [ K[Y,Z](-d_j) ]
$$ 
is a graded free $S$-module of rank $m$ 
and $N$ a graded submodule of $F$. 
An element $\bff \in F$ can be written uniquely as 
$\bff=\ba_0+\ba_{1}Z+\cdots +\ba_dZ^{d}$ with $\ba_i \in \oplus^s K[Y]$.
Denote by $\In'(\bff)$ the term $\ba_j Z^j$ 
for the minimal $j$ such that $\ba _j \not  = {\bf 0}$. 
Furthermore we define 
${\rm In}'(N)$ to be the graded submodule of $F$ 
generated by the set $\{ {\rm In}' (\bff) \}$, 
where $\bff$ runs over homogeneous forms in $N$. 
Also, put 
$$
\Gr_{(Z)}(V)=V/ZV \oplus ZV/Z^2V \oplus Z^2V/Z^3V \oplus \cdots.  
$$
Then, we have that $\Gr_{(Z)}(V) \cong F/\In'(N)$. 
 
Suppose 
$$
F_1 \stackrel{\phi} {\ra} F_0 \ra V \ra 0
$$ 
is a finite presentation of $V$.   
Let $\Delta_{\rm \tiny MAX}(V)$ be the ideal of $S$ 
generated by the maximal minors of $\phi$. 
As a Fitting ideal 
it does not depend on the choice of the finite presentation. 
\end{notation and remark}

\begin{proposition} \label{prop}     
With the same notation as above, 
let $g$ be a general linear form of $S$. 
Put $V^{\ast}=\Gr _ {(Z)}(V)$. 
Assume $K$ is an infinite field. 
Then 
$$
\dim V/gV \leq \dim V^{\ast}/gV^{\ast}. 
$$
\end{proposition}

In the proof of this proposition, 
we use the following two lemmas.

\begin{lemma} \label{lemma-1}               
With the same notation as above, 
let $\ell$ be any linear form of $S$.    
\begin{itemize}
\item[(1)]
$\dim V/\ell V = \dim S/\Delta _{\rm \scriptsize MAX}(V) +(\ell)$. 
\item[(2)]
$\dim V^\ast/\ell V^\ast = \dim S/\Delta _{\rm \scriptsize MAX}(V^\ast) +(\ell)$. 
\end{itemize}
\end{lemma}

\begin{proof} 
(1) 
Obviously $S/\ell S$ is isomorphic to 
the polynomial ring $K[W]$ in one variable.  
Since  $V/\ell V$ is a graded module over a PID $K[W]$, 
we have 
$V/\ell V \cong \oplus_{i=1}^n K[W]/(W^{e_i})$ 
for some positive integers $e_i$ $(1\leq i\leq n)$.
So $\dim V/\ell V$ is equal to $\sum_{i=1}^n e_i$. 
Now, we take a finite presentation of $V/\ell V$, 
$$
\oplus^n K[W] \stackrel{\phi} {\ra} \oplus^n K[W] \ra V/\ell V \ra 0, 
$$
where 
$$
\phi=
\left[
\begin{array}{ccccccc}
W^{e_1}     &         &        & \bigzerou  \\
            & W^{e_2} &        &            \\
            &         & \ddots &            \\
\bigzerol   &         &        & W^{e_n}    \\     
\end{array}
\right]. 
$$
Since $\Delta_{{\rm \scriptsize MAX}}(V/\ell V)=(W^{\sum e_i})$, 
it follows that 
$\dim K[W]/\Delta_{{\rm \scriptsize MAX}}(V/\ell V)=\sum_{i=1}^n e_i$. 
Hence we obtain that 
$\dim V/\ell V=\dim K[W]/\Delta_{{\rm \scriptsize MAX}}(V/\ell V)$. 
So it suffices to show the following  
\medskip

\noindent
{\bf Claim}.   
$K[W]/\Delta_{{\rm \scriptsize MAX}}(V/\ell V) 
\cong  K[Y,Z]/\Delta_{{\rm \scriptsize MAX}}(V)+(\ell).$ 
\medskip

\noindent
{\em Proof of Claim}. 
Suppose 
$$
\oplus^t K[Y,Z] \stackrel{[f_{ij}]} {\ra} \oplus^s K[Y,Z] \ra V \ra 0 
$$
is a finite presentation of $V$, where $f_{ij} \in K[Y,Z]$. 
Taking the tensor product $\otimes_{K[Y,Z]} K[Y,Z]/(\ell)$, 
we get a finite presentation of $V/\ell V $, 
$$
\oplus^t K[W] \stackrel{[f_{ij}']} {\ra} \oplus^s K[W] \ra V/\ell V \ra 0. 
$$
Hence, 
noting that $f_{ij}'$ is the image of $f_{ij}$ in $K[W]$, 
we obtain the desired isomorphism.

We obtain the equality (2) similarly. 
\end{proof}

\begin{lemma} \label{lemma-2}   
Let $I \subset S=K[Y,Z]$ be a homogeneous ideal 
and let $g \in S_1$ be a general linear form. 
Then $\dim S/I+(g)$ is equal to 
the least degree of the monomials that appear in the generators of $I$.  
\end{lemma}

\begin{proof}
Left to the reader. 
\end{proof}

\begin{proof}[Proof of Proposition~\ref{prop}] 
We use Notation~\ref{notation and remark 3-2}. 
For any element $\bff \in F$, 
put $\bff^{\ast}=\In'(\bff)$. 
We can choose generators $\bff_1, \ldots, \bff_t$ for $N$ 
so that $\In'(N)=S\bff_1^{\ast} + \cdots + S\bff_t^{\ast}$. 
We may regard the set of generators $\phi=(\bff_1, \cdots, \bff_t)$ 
as a matrix giving a finite presentation for $V$. 
That is, if we 
let $\bff_i=(f_{i1},\ldots,f_{is})$ 
with homogeneous forms $f_{ij} \in K[Y,Z]$, 
then 
$$
\oplus^t K[Y,Z] \stackrel{\phi} {\ra} \oplus^s K[Y,Z] \ra V \ra 0 
$$
is a finite presentation of $V$, 
where 
$$
\phi=
\left[
\begin{array}{ccccccc}
f_{11}  & \cdots & f_{1s}  \\
\vdots  &        & \vdots  \\
f_{t1}  & \cdots & f_{ts}  \\     
\end{array}
\right]. 
$$
Each $\bff_i$ can be written uniquely as 
$$
\bff_i=(f_{i1}^{(d_i)},\ldots,f_{is}^{(d_i)})Z^{d_i}+
       (f_{i1}^{(d_i+1)},\ldots,f_{is}^{(d_i+1)})Z^{d_i+1}+\cdots, 
$$
where $f_{ij}^{(k)}\in K[Y]$ and 
$$
\bff_i^\ast=(f_{i1}^{(d_i)},\ldots,f_{is}^{(d_i)})Z^{d_i}
         =(f_{i1}^{(d_i)}Z^{d_i},\ldots,f_{is}^{(d_i)}Z^{d_i}). 
$$
Similary, 
we may regard the set of generators $\phi'=(\bff_1^{\ast}, \cdots, \bff_t^{\ast})$ 
as a matrix giving a finite presentation for $V^{\ast}$, 
$$
\oplus^t K[Y,Z] \stackrel{\phi'} {\ra} \oplus^s K[Y,Z] \ra V^{\ast} \ra 0,  
$$
where 
$$
\phi^\prime=
\left[
\begin{array}{ccccccc}
f_{11}^{(d_1)}Z^{d_1}  & \cdots & f_{1s}^{(d_1)}Z^{d_1}  \\
\vdots  &        & \vdots  \\ 
f_{t1}^{(d_t)}Z^{d_t} & \cdots &  f_{ts}^{(d_t)}Z^{d_t}  \\     
\end{array}
\right]. 
$$
Let $\Delta$ be a maximal minor of $\phi$ 
and let $\Delta^\ast$ be the corresponding maximal minor of $\phi^\prime$. 
Note both $\Delta$ and $\Delta^\ast$ are homogeneous forms. 
Obviously, $\Delta^\ast \in K[Y]Z^d$ 
where $d=d_1+\cdots+d_t$.  
Noting that 
$$
f_{ij}=f_{ij}^{(d_i)}Z^{d_i}+f_{ij}^{(d_{i}+1)}Z^{d_i+1}+\cdots, 
$$
it follows that 
$\Delta$ can be written uniquely as $\Delta=\Delta^\ast+\Delta^\prime$, 
where $\Delta^\prime\in\oplus_{m>d}K[Y]Z^m$. 
Hence, the set of monomials that appear in $\Delta^\ast$ 
is contained in that of $\Delta$. 
Thus, the set of monomials that appear in 
the generators of $\Delta _{\rm \scriptsize MAX}(V^{\ast})$ 
is contained in that of $\Delta _{\rm \scriptsize MAX}(V)$. 
>From Lemma~\ref{lemma-1}, it follows that 
$\dim V/g V = \dim S/\Delta _{\rm \scriptsize MAX}(V) +(g)$ 
and $\dim V^{\ast}/g V^{\ast} = 
\dim S/\Delta _{\rm \scriptsize MAX}(V^{\ast}) +(g)$. 
Thus we obtain from Lemma~\ref{lemma-2} that 
$\dim V/gV \leq \dim V^{\ast}/gV^{\ast}$. 
\end{proof}

\begin{lemma} \label{lemma-3}  
Let $A=\oplus_{i=0}^c A_i$ be a standard graded Artinian algebra 
with the WLP (resp. SLP) 
and let $y, z \in A_1$ be two linear forms of $A$. 
If $y$ is a weak Lefschetz element (resp. strong Lefschetz element) of $A$, 
then so is $y+\lambda z$ for a general element $\lambda\in K$.  
\end{lemma}

\begin{proof}  
This follows from the semi-continuity of rank. 
\end{proof}

\begin{notation and remark} \label{notation and remark 3-3}  
Let $A=\oplus_{i=0}^c A_i$ be a standard graded Artinian $K$-algebra 
and let $z$ be a linear form of $A$. 
Since $A$ is Artinian, 
the linear map $\times z: A \ra A$ is nilpotent. 
Hence the Jordan canonical form $J$ of $\times z$ is the matrix as follows. 
$$
J=
\left[
\begin{array}{ccccccc}
J(0,n_1) &     &        & \bigzerou  \\
    & J(0,n_2) &        &            \\
    &     & \ddots &            \\
\bigzerol   &      &        & J(0,n_r)   \\     
\end{array}
\right],  
$$
where $n_1\geq n_2\geq \cdots \geq n_r$ 
and $J(0,m)$ is the Jordan block of size $m\times m$, that is, 
$$
J(0,m)=
\left[
\begin{array}{ccccccc}
0 &  1   &        & \bigzerou  \\
    & 0 &  \ddots      &            \\
    &     & \ddots &    1        \\
\bigzerol   &      &        & 0   \\     
\end{array}
\right].  
$$
Here we note 
$n_1=\Min\{\ i \mid z^i=0 \}$ and $r=\dim_K A/(z)$. 
Now we can take $r$ elements 
$a_i\in A\backslash(z)$ $(i=1,2,\ldots,r)$ such that 
the set 
\begin{equation} \label{EQ:3-7-1}
\cup_{i=1}^r\{a_i,a_iz,a_iz^2,\ldots,a_iz^{n_i-1}\}
\end{equation}
is a basis for $A$ as a vector space, that is, 
the matrix of $\times z:A\ra A$ for the basis above coincides with $J$. 
Then it is easy to check that the set 
$$
\cup_{i=1}^r\{a_i^\ast,a_i^\ast z^\ast,a_i^\ast(z^\ast)^2,\ldots,
a_i^\ast(z^\ast)^{n_i-1}\}
$$
is a basis for $\Gr_{(z)}(A)$. 
Hence the Jordan canonical form of $\times z^{\ast}:\Gr_{(z)}(A)\ra \Gr_{(z)}(A)$ is 
the same as that of $\times z:A\ra A$. 
\end{notation and remark}

\begin{proof}
[Proof of Theorem~\ref{main-th1}]
($\Rightarrow$)\ 
Let $z$ be a  Lefschetz element of $A$. 
Recall that  $A$ and $\Gr_{(z)}(A)$ have the same Hilbert function. 
Then Remark~\ref{notation and remark 3-3} and Lemma~\ref{lemma3} shows that 
$z^{\ast}$ is a Lefschetz element of $\Gr_{(z)}(A)$.

($\Leftarrow$)\ 
We may assume that 
$z$ is the image of $x_n$ in $A=R/I$. 
Then we have that $\Gr_{(z)}(A) \cong R/\In'(I)$. 
\medskip

\ni
{\em Step 1: }  
First we prove that if $\Gr_{(z)}(A)$ has the WLP then so does $A$. 
Let $y^{\ast}$ be a general linear form of $\Gr_{(z)}(A)$, 
let $G \in R_1$ be a preimage of $y^{\ast}$ and 
let $y$ be the image of $G$ in $A=R/I$. 
Let $S=K[Y,Z]$ be the polynomial ring in two variables. 
Define the algebra homomorphism $S \ra \End(A)$  
by $Y \mapsto \times y$ and $Z \mapsto \times z$. 
Then we may regard $A$ as a graded Artinian $S$-module. 
Then, from Proposition~\ref{prop}, 
it follows that 
\begin{equation}\label{EQ}
\dim A/(y+\lambda z)A \leq 
\dim \Gr_{(z)}(A)/(y^{\ast}+\lambda z^{\ast}) \Gr_{(z)}(A)
\end{equation}
for a general element $\lambda \in K$, 
where $z^\ast=z+(z^2)$ in $\Gr_{(z)}(A)$. 

On the other hand since $\Gr_{(z)}(A)$ has the WLP by assumption, 
it follows from Lemmas~\ref{lemma0} (3) and \ref{lemma-3} that 
$$
\dim \Gr_{(z)}(A)/(y^{\ast}+\lambda z^{\ast}) \Gr_{(z)}(A) = 
{\rm Sperner}(\Gr_{(z)}(A)) 
$$ 
for a general element $\lambda\in K$. 
Also it follows that ${\rm Sperner}(A) \leq \dim A/(y+\lambda z)A$ 
from Lemma~\ref{lemma0} (1). 
Furthermore, since $A$ and $\Gr_{(z)}(A)$ have the same Hilbert function, 
the Sperner number of $A$ is equal to that of $\Gr_{(z)}(A)$. 
Hence it follows from the inequality (\ref{EQ}) that
$\dim A/(y+\lambda z)A$ is equal to the Sperner number of $A$. 
Thus, from Lemma~\ref{lemma0} (3), 
we have shown that $A$ has the WLP. 
\medskip

\ni
{\em Step 2:} 
It still remains to show that $A$  has the SLP 
assuming that $\Gr_{(z)}(A)$ has the SLP. 
Let $u$ be a new variable and let $\widetilde{A}=A[u]/(u^{\alpha})$, 
where $\alpha$ is any positive integer. 
Then since we have 
$$\Gr_{(z)}(\widetilde{A})\cong \Gr_{(z)}(A) \otimes _K K[u]/(u^{\alpha})$$
and since the SLP is preserved by tensor product,   
it follows that $\Gr_{(z)}(\widetilde{A})$ has the SLP. 
By Step 1, 
this implies that 
$\widetilde{A}$ has the WLP for all $\alpha > 0$.
Hence the SLP of $A$ follows by Lemma~\ref{lemma4}. 
\end{proof}

\begin{remark} 
Let ${\rm In}(I)$ be the initial ideal of 
a homogeneous ideal $I$ of $R=K[x_1,\ldots,x_n]$ 
with respect to 
the reverse lexicographic order with $x_1>\cdots>x_n$. 
Using a result of \cite{aC03}, 
it can be proved that if $R/{\rm In}(I)$ has the the SLP 
then so does $R/I$ (\cite{aW03}, Proposition 2.9).  
The following example shows that 
Wiebe's result does not imply Theorem~\ref{main-th1} of this paper. 

Let $R=K[x,y,z]$ and $I=(x^2, (x+y)^2, (x+y+z)^2)$. 
Put $A=R/I$. 
Then  
$$
\Gr_{z}(A) \cong R/(x^2, 2xy+y^2, xz+yz, y^3, y^2z, z^3)
$$ 
has the SLP, but 
$$
R/{\rm In}(I)=R/(x^2, xy, xz, y^3, y^2z,z^3)
$$ 
does not have the SLP 
(this can be checked, for instance, with the computer program Macaulay). 
\end{remark}

%
%
%
%

\section{A basic property of central simple modules}

\begin{notation} \label{Jordan sequence}   
Let $A=\oplus_{i=0}^c A_i$ be a standard graded Artinian $K$-algebra 
and let $z$ be a linear form of $A$. 
With Notation~\ref{notation and remark 3-3}, 
let $(f_1,f_2,\ldots,f_s)$ be the finest subsequence of $(n_1,n_2,\ldots,n_r)$ 
such that $f_1>f_2>\cdots>f_s$. 
Then we rewrite the same sequence 
$(n_1, \ldots, n_r)$ as
$$
(n_1, \ldots, n_r)=(\underbrace{f_1 , \ldots f_1}_{m_1}, 
\underbrace{f_2, \ldots, f_2}_{m_2}, \ldots , 
\underbrace{f_s, \ldots, f_s}_{m_s}). 
$$
We call this the {\em Jordan-block-size} 
of the endomorphism $\times z \in \End(A)$. 
\end{notation}

\begin{definition}   
Let the notation be as above. 
We call the graded $A$-module 
$$
U_i=\displaystyle\frac{(0:z^{f_i})+(z)}{(0:z^{f_{i+1}})+(z)} 
$$ 
the {\em $i$-th central simple module} of $(A,z)$, 
where $1\leq i\leq s$ and $f_{s+1}=0$.  
Note these are defined for a pair of the algebra $A$ 
and a linear form $z\in A_1$. 
\end{definition}

\begin{remark}   
By the definition, 
it is easy to see that 
the modules $U_1,U_2,\ldots,U_s$ are non-zero terms of 
the successive quotients of the descending chain of ideals 
$$
A=(0:z^{f_1})+(z) \supset (0:z^{f_1-1})+(z) \supset \cdots 
\supset (0:z)+(z) \supset (z). 
$$
\end{remark}

The Hilbert function of a graded vector space $V=\oplus_{i=a}^b V_i$ 
is the map $i \mapsto \dim V_i$.  
If $V$ has finite dimension, 
then its Hilbert series is the polynomial 
$$
h_V(q) = \sum_{i=a}^b(\dim V_i)q^i. 
$$
Let $h(q)$ be a polynomial with coefficients of integers. 
We say that $h(q)$ is {\em symmetric} 
if $h(q) = q^d h(q^{-1})$ for some integer $d$. 
Then we call the half integer $d/2$ 
the {\em reflecting degree} of the symmetric polynomial $h(q)$. 

The following two lemmas are easy, 
so we omit the proofs.

\begin{lemma} \label{lemma5-2}   
Let $h(q)$ be a polynomial with coefficients of integers. 
The following conditions are equivalent. 
\begin{itemize}
\item[(i)]
$h(q)$ is symmetric. 
\item[(ii)]
$h(q)(1+q+q^2+\cdots+q^k)$ is symmetric 
for all $k=0,1,2,\ldots$. 
\item[(iii)]
$h(q)(1+q+q^2+\cdots+q^k)$ is symmetric 
for some $k\geq 0$. 
\end{itemize}
\end{lemma}

\begin{lemma} \label{lemma5-3}   
Let $h_1(q), h_2(q)$ and $h_3(q)$ 
be polynomials with coefficients of integers. 
Suppose that $h_1(q)=h_2(q)-h_3(q)$, 
$h_2(q)$ is symmetric 
and $h_3(q)$ is symmetric  
with the same reflecting degree as that of $h_2(q)$. 
Then $h_1(q)$ is also symmetric  
with the same reflecting degree as that of $h_2(q)$.
\end{lemma}

For a standard graded Artinian $K$-algebra $A=\oplus_{i\geq 0}A_i$ 
and a homogeneous ideal $J=\oplus_{i\geq 0}J_i$ of $A$, define 
$$
\sigma(A) = 1+{\rm Max}\{i \mid A_i\neq (0) \}  
\ \ \ \mbox{and}\ \ \ 
\alpha(J) = {\rm Min}\{i \mid J_i\neq (0) \}. 
$$

\begin{proposition}  \label{important_property} 
With the notation~\ref{Jordan sequence}, 
suppose that $A$ is an Artinian Gorenstein $K$-algebra 
and let $z$ be a linear form of $A$. 
Let $U_1,\ldots,U_s$ be the central simple modules of $(A,z)$, so  
$$
U_i=\displaystyle\frac{(0:z^{f_i})+(z)}{(0:z^{f_{i+1}})+(z)}  
$$
for all $1\leq i\leq s$. 
Put $\ti{U_i}=U_i\otimes_K K[t]/(t^{f_i})$ for $1\leq i\leq s$. 
Then we have the following. 
\begin{itemize}
\item[(1)]
$h_{\ti{U_i}}(q) = 
h_{U_i}(q)(1+q+q^2+\cdots+q^{f_i-1})$.
\item[(2)]
$h_A(q) = \sum_{i=1}^s h_{\ti{U_i}}(q)$. 
\item[(3)]
$h_{U_i}(q)$ is symmetric for all $i=1,2,\ldots,s$. 
\item[(4)]
$h_{\ti{U_i}}(q)$ is symmetric for all $i=1,2,\ldots,s$ 
with the same reflecting degree as that of $h_A(q)$. 
\item[(5)] 
If all $h_{\ti{U_i}}(q)$ are unimodal, 
then the Sperner number of $A$ is 
the sum of the Sperner numbers of $\ti{U_i}$. 
\end{itemize}
\end{proposition}

\begin{proof}
Since $\ti{U_i}=U_i\otimes_K K[t]/(t^{f_i})$ by definition, 
(1) is immediate. 
Let $\ti{U}=\oplus_{i=1}^s\ti{U_i}$. 
We use Notation~\ref{notation and remark 3-3}. 
Consider the following set, 
\begin{equation} \label{EQ:4-6-1}
\cup_{i=1}^r\{\overline{a_i}, \overline{a_i}\otimes \overline{t}, 
\overline{a_i}\otimes \overline{t}^2,\ldots,\overline{a_i}\otimes \overline{t}^{n_i-1}\}. 
\end{equation}
Then, since the sets (\ref{EQ:3-7-1}) and (\ref{EQ:4-6-1}) 
are bases for $A$ and $\ti{U}$ as graded vector spaces respectively, 
we have that $A$ and $\ti{U}$ have the same Hilbert function. 
Hence, noting that $\ti{U}=\oplus_{i=1}^s\ti{U_i}$, 
we get the equality (2). 
(3) follows from (1) and (4) by Lemma~\ref{lemma5-2}. 
To prove (4) we induct on the number $s$ of central simple modules. 
We use the following well known facts 
on Artinian Gorenstein $K$-algebras. 
\begin{itemize}
\item[1.]
Any Artinian Gorenstein $K$-algebra has a symmetric Hilbert function. 
\item[2.]
If $A$ is an Artinian Gorenstein $K$-algebra, 
then $A/0:a$ is an Artinian Gorenstein $K$-algebra 
for any element $a\in A\backslash\{0\}$ (\cite{jW73}, Lemma 4).   
\end{itemize}
If $s=1$, 
Claim (4) is clear by Fact~1. 
Suppose for the moment $s>1$. 
Let $\overline{A}=A/0:z^{f_s}$ 
and $\overline{z}$ the image of $z$ in $\overline{A}$. 
Then it follows that the Jordan-block-size of 
$\times \overline{z}:\overline{A}\ra\overline{A}$ 
is given by 
$$
(\underbrace{f_1-f_s , \ldots f_1-f_s}_{m_1}, 
\underbrace{f_2-f_s, \ldots, f_2-f_s}_{m_2}, \ldots , 
\underbrace{f_{s-1}-f_s, \ldots, f_{s-1}-f_s}_{m_{s-1}}). 
$$ 
Here we note that 
the image of the ideal $0:z^k$ in $\overline{A}$ 
is equal to 
$(0:\overline{z}^{k-f_s})$ for all $k> f_s$, 
and consider the following descending chain in $\overline{A}$
$$
\overline{A}=(0:\overline{z}^{f_1-f_s})+(\overline{z}) 
\supset (0:\overline{z}^{f_1-f_s-1})+(\overline{z}) 
\supset \cdots 
\supset (0:\overline{z})+(\overline{z}) 
\supset (\overline{z}). 
$$
It is easy to see that 
the central simple modules of $\overline{A}$ 
are the non-zero terms of 
the successive quotients of the chain above. 
Hence they are isomorphic to $U_1,\ldots,U_{s-1}$, respectively. 
Let $\ti{W_i}=U_i\otimes K[t]/(t^{f_i-f_s})$ 
for all $i=1,2,\ldots,s-1$. 
By induction hypothesis, 
we have that 
$h_{\ti{W_i}}(q)$ is symmetric for all $i=1,2,\ldots,s-1$ 
with the same reflecting degree as that of $h_{\overline{A}}(q)$, 
$$
\mbox{(the reflecting degree of $h_{\ti{W_i}}(q)$)} 
= (\sigma(\overline{A})-1)/2. 
$$
>From Theorem 3 (c) in \cite{DGO}, 
it follows that 
$$
\sigma(\overline{A})=\sigma(A)-\alpha((z^{f_s}))=\sigma(A)-f_s.  
$$
Also, since 
$$
h_{\ti{W_i}}(q) = 
h_{U_i}(q)(1+q+q^2+\cdots+q^{f_i-f_s-1}) 
$$
for all $i=1,2,\ldots,s-1$, 
it follows from (1) that 
$$
\mbox{(the reflecting degree of $h_{\ti{U_i}}(q)$)} 
= \mbox{(the reflecting degree of $h_{\ti{W_i}}(q)$)} +f_s/2. 
$$
Hence
$$
\begin{array}{rcl}
\mbox{(the reflecting degree of $h_{\ti{U_i}}(q)$)} 
  & = & (\sigma(A)-1)/2 \\[1ex]
  & = & 
\mbox{(the reflecting degree of $h_{A}(q)$)} 
\end{array}
$$
for all $i=1,2,\ldots,s-1$. 
Furthermore, since 
$$
h_{\ti{U_s}}(q) = h_A(q) - \sum_{i=1}^{s-1} h_{\ti{U_i}}(q), 
$$
it follows from Lemma~\ref{lemma5-3} that 
$h_{\ti{U_s}}(q)$ is also symmetric  
with the same reflecting degree as that of $h_A(q)$. 

Finally (5) follows from (4) and (2). 
\end{proof}

\begin{definition}   
Let $A=\oplus_{i\geq 0} A_i$ be any graded $K$-algebra. 
Suppose that 
$$
V=\oplus_{i=a}^b V_i
$$
is a graded finite $A$-module with $V_a\neq(0)$ and $V_b\neq(0)$. 
\begin{itemize}
\item[(i)]
The $A$-module $V$ has the {\em WLP} 
if there is a linear form $g\in A_1$ such that 
the multiplication $\times g: V_i\ra V_{i+1}$ is 
either injective or surjective for all $i=a,a+1,\ldots,b-1$. 
\item[(ii)]
The $A$-module $V$ has the {\em SLP} 
if there is a linear form $g\in A_1$ such that 
the multiplication $\times g^{b-a-2i}: V_{a+i}\ra V_{b-i}$ 
is bijective for all $i=0,1,\ldots,[(b-a)/2]$. 
\end{itemize} 
\end{definition}

The following is immediate by Proposition~\ref{important_property}, so we omit the proof.

\begin{corollary} \label{cor4.7} 
In addition to the assumption of Proposition~\ref{important_property}, 
assume furthermore that all central simple modules of $(A,z)$ have the SLP. 
Then $\ti{U}=\oplus_{i=1}^s\ti{U_i}$ has the SLP. 
\end{corollary}

%
%
%
%

\section{Proof of Main Theorem~\ref{main-th2}}

\begin{proof}[Proof of Theorem~\ref{main-th2}] 
(i)$\Rightarrow$(ii): 
Assume that $A$ has the SLP and $z$ is a strong Lefschetz element. 
Put 
$\overline{A}_i =A_i/zA_{i-1}$.  
Then we may write 
$$A/(z)= \oplus _{i=0}^ {c'}\overline{A_i}$$
where $c'$ is the largest integer such that 
$(A/(z))_{c'} \not = 0$. 
Then, noting that $A$ has the SLP, 
one sees easily that $(A,z)$ has  
$c'+1$ central simple modules $U_1, \cdots, U_{c'+1}$ 
and that $U_i= \overline{A}_{i-1}$.  
In particular this shows that $U_i$ has only one non-trivial 
graded piece concentrated at the degree $i-1$, hence 
$U_i$ has the SLP for trivial reasons. 
\medskip

\ni
(ii)$\Rightarrow$(i): 
We use Notation~\ref{Jordan sequence}. 
We divide the proof into three steps. 
\medskip

\ni
{\em Step1}: 
Since the set of strong Lefschetz elements is a Zariski open subset of $A_1$, 
there exists a linear form $g \in A$ such that 
$g$ is a strong Lefschetz element of $U_i$ for every $i=1,2,\ldots,s$.   
Let $m_i=\dim U_i$ $(1\leq i\leq s)$. 
Now we take a basis of $U_i$ as a vector space for all $1\leq i\leq s$, 
$$
\{ \ol{e_{i1}}, \ol{e_{i2}}, \ldots, \ol{e_{im_{i}}} \}, 
$$
where $e_{ij}\in (0:z^{f_i})+(z)$.  
Then the set 
$$
\{ \ol{e_{ij}}\otimes \overline{t}^k \mid 1\leq j\leq m_i, \ 0\leq k\leq f_i-1\}
$$
is a basis of $\ti{U_i}=U_i\otimes_K K[t]/(t^{f_i})$ 
for all $1\leq i\leq s$. 
Put $\ti{U}=\oplus_{i=1}^s\ti{U_i}$. 
The set 
\begin{equation} \label{basis1}
\cup_{i=1}^s\{ \ol{e_{ij}}\otimes \overline{t}^k \mid 1\leq j\leq m_i, \ 0\leq k\leq f_i-1\}
\end{equation}
is a basis of $\ti{U}$. 
We may consider $\ti{U}$ as a graded $A\otimes K[t]$-module. 
Here we calculate a matrix of the multiplication 
$\times (g\otimes 1 + 1\otimes t) : \ti{U} \ra \ti{U}$ 
which is a linear map.  
Let $P_i$ be the square matrix of 
$\times (g\otimes 1 + 1\otimes t) : \ti{U_i}\ra\ti{U_i}$ 
for the basis above. 
Since $(g\otimes 1 + 1\otimes t)\ti{U_i}\subset\ti{U_i}$, 
a matrix of $\times (g\otimes 1 + 1\otimes t) :\ti{U} \ra \ti{U}$ 
is of the following form,  
$$
P=
\left[
\begin{array}{ccccccc}
P_1 &     &        & \bigzerou  \\
    & P_2 &        &            \\
    &     & \ddots &            \\
\bigzerol   &      &        & P_s   \\     
\end{array}
\right]. 
$$
Hence it follows that 
$$
P^h=
\left[
\begin{array}{ccccccc}
P_1^h &     &        & \bigzerou  \\
    & P_2^h &        &            \\
    &     & \ddots &            \\
\bigzerol   &      &        & P_s^h   \\     
\end{array}
\right]
$$
for all $h$. 
\medskip

\ni
{\em Step 2}: 
Let $g^\ast$ be the initial form of $g$ in $\Gr_{(z)}(A)$. 
We calculate a matrix of the multiplication 
$\times g^\ast: \Gr_{(z)}(A)\ra \Gr_{(z)}(A)$. 
First we note that the set 
\begin{equation} \label{basis2}
\cup_{i=1}^s \{ (e_{ij}^\ast)(z^\ast)^k \mid 1\leq j\leq m_i, 
\ 0\leq k\leq f_i-1 \}
\end{equation}
is a basis of $\Gr_{(z)}(A)$. 
Let $V_i$ be the subspace of $A/(z)$ in $\Gr_{(z)}(A)$ 
which is generated by $\{e_{i1}^\ast, e_{i2}^\ast, \ldots, e_{im_i}^\ast \}$ 
for all $i=1,2,\ldots,s$. 
Furthermore, 
let $V_i^\ast$ be a subspace of $\Gr_{(z)}(A)$ which is generated by 
$$
\{ (e_{ij}^\ast)(z^\ast)^k \mid 1\leq j\leq m_i, \ 0\leq k\leq f_i-1 \}
$$ 
for all $i=1,2,\ldots,s$. 
Then, since $g^\ast V_i\subset\oplus_{j=i}^s V_i$ for all $i=1,2,\ldots,s$, 
we have that $g^\ast V_i^\ast\subset\oplus_{j=i}^s V_j^\ast$.  
Hence, a matrix of 
$\times (g^\ast + z^\ast) : \Gr_{(z)}(A)\ra \Gr_{(z)}(A)$ is 
of the following form 
$$
Q=
\left[
\begin{array}{ccccccc}
P_1 &     &        & \bigast  \\
    & P_2 &        &            \\
    &     & \ddots &            \\
\bigzerol   &      &        & P_s   \\     
\end{array}
\right]. 
$$
Thus, we obtain that 
$$
Q^h=
\left[
\begin{array}{ccccccc}
P_1^h &     &        & \bigast  \\
    & P_2^h &        &            \\
    &     & \ddots &            \\
\bigzerol   &      &        & P_s^h   \\     
\end{array}
\right]
$$
for all $h$. 
\medskip

\ni
{\em Step 3}: 
First note the following. 
\begin{itemize}
\item[(1)]
$\ti{U}$ and $\Gr_{(z)}(A)$ have the same Hilbert function. 
To see this note that  $\deg\ol{e_{ij}}=\deg e_{ij}^\ast$ for all $i$ and $j$ 
and the sets (\ref{basis1}) and (\ref{basis2}) 
are bases for $\ti{U}$ and $\Gr_{(z)}(A)$ respectively. 
\item[(2)] 
$\ti{U}$ has the SLP. This follows from Corollary~\ref{cor4.7}. 
\end{itemize}

Hence we have from arguments analogous to 
Lemma \ref{lemma2} (2) and (4) that 
if $\rank \ P^h \leq \rank \ Q^h$ for all $h$,  
then $\Gr_{(z)}(A)$ has the SLP, which further enables us to conclude that 
$A$ has the SLP, thanks to Theorem~\ref{main-th1}. 
So it suffices to show that $\rank \ P^h \leq \rank \ Q^h$ for all $h$. 
Let $P_{i h}^\prime$ be a square submatrix of $P_i^h$ 
such that $P_{i h}^\prime$ is of full rank  and 
$\rank \ P_i^h = \rank \ P_{i h}^\prime$ 
for all $i$ and $h$. 
Then we have  
$$
\begin{array}{rcl}
\rank \ P^h & = & \sum_{i=1}^s \rank \ P_{i}^h \\[1ex]
& = & \sum_{i=1}^s \rank \ P_{i h}^\prime \\[1ex]
& = & \rank \ P_h^\prime,  
\end{array}
$$ 
where 
$$
P_h^\prime=
\left[
\begin{array}{ccccccc}
P_{1 h}^\prime &     &        & \bigzerou  \\
    & P_{2 h}^\prime &        &            \\
    &     & \ddots &            \\
\bigzerol   &      &        & P_{s h}^\prime   \\     
\end{array}
\right]. 
$$
Let $Q_h^\prime$ be the square submatrix of $Q^h$ 
consisting of the same rows and columns as $P_h^\prime$, 
so that  
$Q_h^\prime$ is of  the following form 
$$
Q_h^\prime=
\left[
\begin{array}{ccccccc}
P_{1 h}^\prime &     &        & \bigast  \\
    & P_{2 h}^\prime &        &            \\
    &     & \ddots &            \\
\bigzerol   &      &        & P_{s h}^\prime   \\     
\end{array}
\right]. 
$$
Then, since $\det Q_h^\prime =\prod_{i=1}^s \det P_{i h}^\prime \neq 0$, 
it follows that $\rank \ Q_h^\prime = \rank \ P_h^\prime = \rank \ P^h$. 
This means that $\rank \ P^h \leq \rank \ Q^h$ for all $h$. 
\end{proof}

%
%
%
%

\section{Some examples of complete intersections with the SLP}

\begin{lemma} \label{lemma6-1}   
Let $R=K[x_1, \ldots, x_n]$ be the polynomial ring over a field $K$ 
and let $J$ be a homogeneous ideal of $R$ 
such that $R/J$ is a one dimensional Cohen-Macaulay $K$-algebra. 
Let $g$ be a linear form of $R$ which is not a zero divisor on $R/J$. 
Let $d$ be a positive integer, 
and put $I=(J,g^d)$ and $A=R/I$. 
Let $z$ be the image of $g$ in $A$. 
Then $(A,z)$ has only one central simple module 
which is isomorphic to $A/(z)$. 
\end{lemma}

\begin{proof}
It is easy to see  that 
$I:g^j = (J, g^{d-j})$ 
for all $j=0,1,\ldots,d$. 
Hence, since $(0:z^j) = (z^{d-j})$ 
for all $j=0,1,\ldots,d$, 
we have 
$(0:z^j)+(z) = (z)$ 
for all $j=0,1,\ldots,d-1$. 
Thus $(A,z)$ has only one central simple module 
$$
\displaystyle\frac{(0:z^d)+(z)}{(0:z^{d-1})+(z)}=\frac{A}{(z)}. 
$$ 

\end{proof}

\begin{example} \label{example1}   
Let $R=K[x_1, \cdots, x_n]$ be the polynomial ring over a field $K$ of characteristic~0. 
We consider a complete intersection ideal as follows, 
$$
I=(f_1, f_2, g_3^{d_3}, g_4^{d_4}, \ldots, g_n^{d_n}), 
$$
where $g_3,\ldots,g_n$ are linear forms. 
Then $A=R/I$ has the SLP. 
\end{example}

\begin{proof}
If $n \leq 2$,  this is proved in  Proposition 4.4 of \cite{tHjMuNjW01}. 
Now we induct on $n$. 
Let $n\geq 3$. 
We may assume that $g_n=x_n$. 
Put $\overline{R}=R/x_nR$. 
Then 
$$
A/x_nA=R/(I,x_n)=\overline{R}/(\overline{f_1}, \overline{f_2}, \overline{g_3}^{d_3}, \ldots, \overline{g_{n-1}}^{d_{n-1}}), 
$$
where $\overline{f_i}$ and $\overline{g_j}$ are the images of $f_i$ and $g_j$ in $\overline{R}$. 
By the assumption of induction, 
we have that $A/x_nA$ has the SLP. 
Hence, our assertion follows 
from Lemma~\ref{lemma6-1} and Theorem~\ref{main-th2}. 
\end{proof}

\begin{notation and remark} \label{notation-symmetric}  
Let $e_i=e_i(x_1, \ldots, x_n)$ be 
the elementary symmetric polynomial of degree $i$ 
in the variables $x_1, \ldots, x_n$, i.e., 
$$
e_i(x_1, \ldots, x_n)=
\sum_{j_1<j_2<\cdots<j_i}x_{j_1}x_{j_2}\cdots x_{j_i} 
$$
for all $i=1,2,\ldots,n$. 
Let $r$ and $s$ be two positive integers. 
Put 
$$
\left\{
\begin{array}{l}
  f_i=e_i(x_1^r, \ldots, x_n^r),  \mbox{ for } i=1, \ldots, n-1,  \\
  f_n=e_n(x_1^s, \ldots, x_n^s).
\end{array}
\right.
$$
It is easy to see that  the ideals  $(e_1,e_2,\ldots,e_n)$ 
and   
$(f_1,f_2,\ldots,f_n)$ are complete intersections. 
\end{notation and remark}

\begin{example} \label{example2}   
With the same notation as above, 
let $R=K[x_1, \ldots, x_n]$ be the polynomial ring 
over a field $K$ of characteristic~0.  
Put  $I=(f_1, \ldots, f_n) \mbox{ and } A=R/I$. 
Suppose $s$ is a multiple of $r$.
Then the complete intersection $A$ has the SLP. 
\end{example}

\begin{proof}
Noting 
\begin{equation}\label{eq-symmetric}
(-1)^{n+1}e_n+(-1)^{n}x_ne_{n-1}+(-1)^{n-1}x_n^2e_{n-2}+
\cdots+(-1)^2x_n^{n-1}e_1=x_n^n, 
\end{equation}
there exist polynomials $P_1,\ldots,P_n \in R$ such that 
$$
x_n^{rn}=P_1f_1+\cdots+P_{n-1}f_{n-1}+P_n x_1^r\cdots x_n^r. 
$$ 
Hence, since 
$$
(x_n^{rn})^{s/r}=(P_1f_1+\cdots+P_{n-1}f_{n-1}+P_n x_1^r\cdots x_n^r)^{s/r},
$$ 
it follows that $x_n^{sn}\in I=(f_1,\ldots,f_{n-1}, f_n)$. 
Thus we obtain that  
$$
(f_1,\ldots,f_{n-1}, x_n^{sn})\subset I. 
$$
On the other hand, noting that 
$R/(f_1,\ldots,f_{n-1}, x_n^{sn})$ and $A=R/I$ are both complete intersections with the same Hilbert function, 
we have 
$$
I=(f_1,\ldots,f_{n-1}, x_n^{sn}). 
$$
Let $z$ be the image of $x_n$ in $A$. 
Then we see that  
$$
A/(z) \cong K[x_1,\ldots,x_{n-1}]/(e_1'(x_1^r,\ldots,x_{n-1}^r),\ldots,e_{n-1}'(x_1^r,\ldots,x_{n-1}^r)), 
$$
where $e_i'(x_1,\ldots,x_{n-1})$ is the elementary symmetric function of degree $i$ 
in the variables $x_1, \cdots, x_{n-1}$ for all $1\leq i\leq n-1$. 
Here, inductively we may assume that $A/(z)$ has the SLP. 
Hence the SLP of $A$ follows from Lemma~\ref{lemma6-1} and Theorem~\ref{main-th2}. 
\end{proof}

\begin{example} \label{example3}   
We keep the notation of Example~\ref{example2}. 
Let $z$ be the image of $x_n$ in $A$. 
Then the complete intersection $A/0:z^k$ has the SLP 
for all $1\leq k <sn$. 
\end{example}

\begin{proof}
>From the proof of Example~\ref{example2}, 
the ideal $I$ is generated by a homogeneous regular sequence 
$\{f_1,\ldots,f_{n-1}, x_n^{sn}\}$. 
Hence  
$$
A/0:z^k \cong R/(f_1,\ldots,f_{n-1}, x_n^{sn-k}). 
$$
Thus, from Lemma~\ref{lemma6-1}, 
we obtain that $A/0:z^k$ has the SLP. 
\end{proof}

\begin{proposition} \label{prop6}   
Let $R=K[x_1,\ldots,x_n]$ be the polynomial ring 
over a field $K$ of characteristic zero, 
and $I$ a complete intersection homogeneous ideal of height $n$. 
Put $A=R/I$ and let $z$ be a linear form of $A$. 
Suppose that $A/0:z^k$ is either a complete intersection 
or $(0)$ for every $k\geq 0$. 
Then all central simple modules of $(A,z)$ are principal.  
In this case all the central simple modules  
are isomorphic to some Artinian Gorenstein algebras with a shift of degrees. 
\end{proposition}

To prove this proposition, 
we need a lemma.

\begin{lemma} \label{lemma6-2}   
Let $I$ be a homogenous complete intersection ideal of 
$R=K[x_1,\ldots,x_n]$ of height $n$. 
Let $\ell \in R$ be a linear form. 
Then the following conditions are equivalent. 

\begin{enumerate}
\item[(1)]
$I+(\ell)$ can be generated by $n$ elements.
\item[(2)]
$I\colon \ell$ can be generated by $n$ elements.
\item[(3)]
A multiple of $\ell$ can be a member of a minimal generating set of $I$.
\end{enumerate}
\end{lemma}

\begin{proof}
\noindent
(1) $\Rightarrow$ (3): 
Let $\{f_1,f_2,\ldots,f_n\}$ be a set of homogeneous generators of $I$. 
The assumption means that 
one of the generators of $(f_1, f_2, \cdots, f_{n}, \ell)$ is superfluous. 
Suppose $\ell$ is superfluous. 
Then, noting that the degree of $\ell$ is one, 
it follows that $\ell$ can be a member of a minimal generating set of $I$. 
Otherwise we have, for example, $f_n=a_1f_1+ \cdots + a_{n-1}f_{n-1}+a\ell$ with some 
elements $a_i, a \in R$.  
Then $a\ell$ can replace $f_n$ in the minimal generating set of $I$. 
\medskip

\noindent
(3) $\Rightarrow$ (2): Straightforward. 
\medskip

\noindent
(2) $\Rightarrow$ (1): 
Let $\{g_1,g_2,\cdots,g_n\}$ be 
a set of homogeneous generators of $I:\ell$. 
Since $I \subset I:\ell$, 
there exists the $n\times n$ matrix $M=[a_{ij}]$ such that
$[f_1,\cdots, f_n]=[g_1,\cdots,g_n]M$, 
where $a_{ij}$ is a homogeneous form. 
Here we recall the following well known equality 
on Koszul complexes of complete intersections 
(for example, see page 983 of \cite{B}), 
$$
\sum_{i=1}^n\deg f_i=(\sum_{i=1}^n\deg g_i)+\deg(\det M). 
$$
Also, from Theorem 3 (c) in \cite{DGO}, 
we have  
$$
\sigma(R/I)=\sigma(R/I:\ell)+\alpha((\overline{\ell}))=\sigma(R/I:\ell)+1,   
$$
where $\overline{\ell}$ is the image of $\ell$ in $R/I$. 
Hence, 
noting $\sigma(R/I)=(\sum_{i=1}^n\deg f_i)-n+1$ and 
$\sigma(R/I:\ell)=(\sum_{i=1}^n\deg g_i)-n+1$, 
we get from the equalities above that 
$$
\deg(\det M)=(\sum\deg f_i)-(\sum\deg g_i)=\sigma(R/I)-\sigma(R/I:\ell)=1.
$$ 
This means that an entry of $M$, say $a_{st}$,   
 is a linear form 
and the co-factor of $a_{st}$ is a non-zero element in $K$, 
that is, $\tilde{a_{st}}=(-1)^{s+t}\det(M')\in K\backslash\{0\}$, 
where $M'$ is the $(n-1)\times (n-1)$ submatrix of $M$ 
deleting the $s$th row and the $t$th column. 
Let $\tilde{M}=[\tilde{a_{ij}}]$ be the ajoint matrix of $M$. 
Then, noting that 
$$
[f_1,f_2,\ldots,f_n]\ ^{t}\tilde{M} 
= [\det(M)g_1,\det(M)g_2,\ldots,\det(M)g_n], 
$$
we have that 
$$
\tilde{a_{s1}}f_1+\cdots+\tilde{a_{st}}f_t+\cdots+\tilde{a_{sn}}f_n 
= \det(M)g_s. 
$$
Hence, since $\tilde{a_{st}}\in K\backslash\{0\}$, 
this means that 
$$
I=(f_1,\ldots,f_{t-1},\det(M)g_s,f_{t+1},\ldots,f_n). 
$$
Thus, 
noting that $I:(I:\ell)=I+(\det(M))$ (for example, see page 983 of \cite{B}), 
we obtain
$$
I+(\ell)=I+(\det M)=(f_1,\ldots,f_{t-1},\det(M),f_{t+1},\ldots,f_n).  
$$
\end{proof}

\begin{proof}[Proof of Proposition~\ref{prop6}] 
(We use the same letter $z$ for the preimage of $z$ in $R$.)
First we would like to show that 
$I:z^k+(z)/I:z^{k-1}+(z)$ can be generated by 
a single element for any integer $k\geq 0$. 
Put $J=I:z^{k-1}$. 
By assumption $J$ and $J:z=I:z^k$ are complete intersections. 
So by Lemma~\ref{lemma6-2} (3), 
there is a set of generators $g_1,\ldots,g_{n-1},g_n$ 
such that $J=(g_1,\ldots,g_{n-1},zg_n)$ 
and $J:z=(g_1,\ldots,g_{n-1},g_n)$. 
Hence the first assertion follows. 
Suppose that $U$ is a principal central simple module. 
Then, with some integer $k>0$, we have 
the exact sequence 
$$
0 \ra U \ra A/0:z^{k-1}+(z) \ra A/0:z^{k}+(z) \ra 0. 
$$
Hence, 
since $U$ is a principal ideal 
of the complete intersection $B=A/0:z^{k-1}+(z)$, 
there exists some $b\in B$ such that $U\cong B/0:b$, 
up to a degree shift. 
Thus it follows from Lemma 4 of \cite{jW73} that 
$U$ is isomorphic to a Gorenstein algebra. 
\end{proof}

\begin{example} \label{example4} 
{\rm
With the notation~\ref{notation-symmetric}, 
let $R=K[x_1, \cdots, x_n]$ be the polynomial ring 
over a field $K$ of characteristic~0. 
Put  $I=(f_1, \cdots, f_n) \mbox{ and } A=R/I$.  
Suppose $s < r$. 
Then $A$ has the SLP. 
In this case, 
taking the image $z$ of $x_n$ in A, 
it follows that 
$(A, z)$ has two central simple modules $U_1 \mbox{ and } U_2$,   
$$
U_1\cong K[x_1, \cdots, x_{n-1}]
/(\overline{f}_1, \cdots, \overline{f}_{n-2}, 
(x_1\cdots x_{n-1})^s), 
$$
$$
U_2\cong K[x_1, \cdots, x_{n-1}]
/(\overline{f}_1, \cdots, \overline{f}_{n-2}, 
(x_1\cdots x_{n-1})^{r-s}),
$$
where $\overline{f}_j$ is the image of $f_j$ 
in $K[x_1, \cdots, x_{n-1}]$. 
}
\end{example}

\begin{proof}
>From Proposition~\ref{prop6}, 
it is sufficient to show that 
$A/0:z^k$ is either a complete intersection or $(0)$ 
for every $k\geq 0$.
First we show that 
\begin{itemize}
\item[(1)]
$I:x_n^k=(f_1,\ldots,f_{n-1},x_1^s\cdots x_{n-1}^sx_n^{s-k})$ 
for all $0\leq k\leq s$, 
\item[(2)]
$I:x_n^k=(f_1,\ldots,f_{n-2},x_n^{(n-1)r-(k-s)},x_1^s\cdots x_{n-1}^s)$ 
for all $s<k<(n-1)r+s$, 
\item[(3)]
$I:x_n^k=R$ for all $k\geq (n-1)r+s$. 
\end{itemize}
(1) is easy. 
So we give a proof of (2). 
Using the equation~(\ref{eq-symmetric})
we get the following relation 
$$
(-1)^{n+1}x_1^r\cdots x_n^r+(-1)^{n}x_n^rf_{n-1}+
(-1)^{n-1}x_n^{2r}f_{n-2}+\cdots+(-1)^2x_n^{(n-1)r}f_1=x_n^{nr}. 
$$
Hence we have 
\begin{equation}\label{eq-symmetric2}
(-1)^{n+1}x_1^r\cdots x_{n-1}^r+(-1)^{n}f_{n-1}+
(-1)^{n-1}x_n^{r}f_{n-2}+\cdots+(-1)^2x_n^{(n-2)r}f_1=x_n^{(n-1)r}, 
\end{equation}
and 
$$
I:x_n^s
=(f_1,\ldots,f_{n-1},x_1^s\cdots x_{n-1}^s)
=(f_1,\ldots,f_{n-2},x_n^{(n-1)r},x_1^s\cdots x_{n-1}^s).
$$
Therefore, for all $s<k\leq(n-1)r+s$, 
it follows that 
$$
I:x_n^k=(I:x_n^s):x_n^{k-s}
=(f_1,\ldots,f_{n-2},x_n^{(n-1)r-(k-s)},x_1^s\cdots x_{n-1}^s). 
$$
(3) is easy. 

Next we calculate the central simple modules of $(A,z)$. 
Using (1) to (3), we have 
$$
(I:x_n^k)+(x_n) = 
\left\{
\begin{array}{ll}
(f_1,\ldots,f_{n-1},x_n)  
& \ k=0,1,\ldots,s-1, \\[1ex]
(f_1,\ldots,f_{n-1},x_1^s\cdots x_{n-1}^s,x_n)  & \ k=s, \\[1ex] 
(f_1,\ldots,f_{n-2},x_1^s\cdots x_{n-1}^s,x_n)  
& \ k=s+1,s+2,\ldots,(n-1)r+s-1, \\[1ex] 
R & \ k=(n-1)r+s,\ldots. 
\end{array}
\right.
$$
Here, noting the equality~(\ref{eq-symmetric2}), 
it follows that 
$$
(f_1,\ldots,f_{n-1},x_1^s\cdots x_{n-1}^s,x_n)
=(f_1,\ldots,f_{n-2},x_1^s\cdots x_{n-1}^s,x_n).  
$$
Hence we have 
$$
(0:z^k)+(z) = 
\left\{
\begin{array}{ll}
(z)  
& \ k=0,1,\ldots,s-1, \\[1ex]
(\overline{x}_1^s\cdots \overline{x}_{n-1}^s,z)  
& \ k=s,s+1,\ldots,(n-1)r+s-1, \\[1ex] 
A & \ k=(n-1)r+s,\ldots, 
\end{array}
\right.
$$
where $\overline{x}_j$ is the image of $x_j$ in $A$. 
Thus we obtain that 
$$
U_1 = A/(\overline{x}_1^s\cdots \overline{x}_{n-1}^s,z) 
\cong K[x_1, \cdots, x_{n-1}]
/(\overline{f}_1, \cdots, \overline{f}_{n-2}, 
(x_1\cdots x_{n-1})^s)
$$
and
$$
\begin{array}{lcl}
U_2 & = & (\overline{x}_1^s\cdots \overline{x}_{n-1}^s,z)/(z) \\[1ex]
 & \cong & \displaystyle\frac{A/(z)}
          {(0:\overline{x}_1^s\cdots\overline{x}_{n-1}^s)+(z)/(z)} \\[3ex]
 & \cong & R/I+(x_n):x_1^s\cdots x_{n-1}^s \\[1ex] 
 & \cong & K[x_1, \cdots, x_{n-1}]
/(\overline{f}_1, \cdots, \overline{f}_{n-2}, 
(x_1\cdots x_{n-1})^{r-s}).
\end{array}
$$
\end{proof}

\begin{example}     
Let $R=K[x,y,z]$ be the polynomial ring in $x,y,z$, 
where $K$ is a field of characteristic zero.  
Let $a>0$ be a positive integer. 
Let $$I=(p_a(x,y,z), p_{a+1}(x,y,z), p_{a+2}(x,y,z)),$$
where $p_d=x^d+y^d+z^d$ is the power sum of degree $d$.
Then $A=R/I$ has the SLP. 
\end{example}
 
\begin{proof}
First we have to show that the ideal $I=(p_a, p_{a+1}, p_{a+2})$ is a complete intersection. 
Using Newton's identity it can be proved that  $p_i \in I$ for all $i > a+2$.  Thus $I \supset (p_a, p_{2a},  p_{3a})$.  
Now it is easy to see that $(p_1, p_2, p_3)$ is a complete intersection. Hence 
$(p_a, p_{2a}, p_{3a})$ is a complete intersection and  so is I. 

Now put $f=(x-z)(y-z), f'=-x-y+2z$.  Note that we have the relations
$$
\left\{\begin{array}{l}
z^af-(xy)p_a+(x+y)p_{a+1}-p_{a+2}=0, 
\\
\\
2z^af'-(h_{a-1}(x,0,z)+h_{a-1}(0,y,z))f+(x+y-z)p_a-p_{a+1}=0.  
\end{array}
\right.
$$
where $h_d(x,y,z)$ denotes the complete symmetric function of degree $d$, i.e.,    
$$
h_d(x,y,z)=\sum_{c_1+c_2+c_3=d}x^{c_1}y^{c_2}z^{c_3}. 
$$
These identities show the equality of ideals 
$$I=(p_a, p_{a+1}, z^af) \ \ \ \mbox{and}\ \ \ I\colon z^a=(f, p_a, p_{a+1})=(f, p_{a}, f'z^a).$$ 
Furthermore it is easy to see that 
$$I\colon z^{2a}=(f',f,p_a)=(f',f,z^a)\ \ \ \mbox{and}\ \ \ I\colon z^{3a}=(1).$$
This shows that $I\colon z^m$ is a complete intersection for every $1\leq m <3a$. 
By  Proposition~\ref{prop6}  it follows 
that all central simple modules for $(A,z)$ are complete intersections of codimension two. 
On the other hand every complete intersection of codimension two has the SLP (\cite{tHjMuNjW01}, Proposition 4.4). 
So by Theorem~\ref{main-th2} we get the assertion. 
\end{proof}

Let $A$ be a complete intersection $K$-algebra with the SLP. 
Then it is not always true that 
all central simple modules of $(A,z)$ have the SLP 
for all linear forms $z$ of $A$. 
Finally we show such an example.

\begin{example}  
Let $A=R/I$, where 
$R=K[u,v,w,x,y,z], I=(y^2, x^2, w^2, v^3, u^3, z^5-zp)$  and  $p=u^2wx+uvwy+v^2xy$.  
Then $(A, \overline{z})$ has three central simple modules, 
$U_1$, $U_2$ and $U_3$. 
The algebra $A$ has the SLP but $U_3$ does not. 
\end{example} 

\begin{proof}  
By a computer algebra system it is possible to compute $\dim A/(z^i)$. 
This gives us the Jordan-block-size of $\times z$ as follows. 
\begin{equation}   \label{grand_block_decomposition}   
(\underbrace{9 , \ldots 9}_{12}, 
\underbrace{5, \ldots, 5}_{48},  
\underbrace{1, \ldots, 1}_{12}). 
\end{equation}
This shows that $(A, \overline{z})$ has 
three central simples modules $U_1, U_2, U_3$ with 
$\dim U_1 = 12$,   $\dim U_2 = 48$,   $\dim U_3 = 12$.  
The Hilbert functions are known to be 
$$
\left\{
\begin{array}{l}
h_{U_1}= 1+5q+5q^2+q^3, \\
h_{U_2}= 7q^2+17q^3+17q^4+7q^5, \\
h_{U_3}= q^4+5q^5+5q^6+q^7. 
\end{array}
\right.
$$
Furthermore 
$$
\left\{
\begin{array}{l}
h_{\ti{U_1}}= (1+5q+5q^2+q^3)(1+q+ \cdots +q^8), \\
h_{\ti{U_2}}= (7q^2+17q^3+17q^4+7q^5)(1+q+ \cdots + q^4), \\
h_{\ti{U_3}}= q^4+5q^5+5q^6+q^7. 
\end{array}
\right.
$$
Notice that 
$$
R/I:z = R/(y^2, x^2, w^2, v^3, u^3, z^4- p). 
$$
Hence 
$$
R/(I:z)+(z) \cong K[u,v,w,x,y]/(y^2, x^2, w^2, u^3, v^3, p).  
$$
This shows that $U_3$ is principal.  In fact as an $A$-module 
$$
U_3 = 0: \overline{z}/ (\overline{z})  \cong \overline{p}+(\overline{z})/ (\overline{z}). 
$$
It is not difficult to see that 
$$
\overline{p}+(\overline{z})/ (\overline{z}) \cong (K[u,v,w,x,y]/{\rm Ann}(F))(-4)
$$
\def\pa{{\tiny \partial}}
where $F=wu^2 + 2xuv + yv^2$ and 
$$
{\rm Ann}(F)=\{f \in K[u,v,w,x,y]| 
f(\frac{\pa}{\pa u}, \frac{\pa}{\pa v}, \frac{\pa}{\pa w}, \frac{\pa}{\pa x}, \frac{\pa}{\pa y})F=0 \}.
$$
It is known that $F$ in Macaulay's inverse system gives a Gorenstein algebra which does  
not have the SLP (\cite{tHjMuNjW01}, Example~4.3). 
That the algebra $A$ has the SLP can be proved by the same method as Example~\ref{example1}.  
\end{proof}
\medskip


\ni
{\bf Acknowledgement}
\medskip

The second author would like to thank Aldo Conca for helpful conversations 
about power sum symmetric functions.
\medskip


\end{document}